\documentclass[12pt]{amsart}

\usepackage{amssymb}
\usepackage{latexsym}
\usepackage{amsmath}

\begin{document}

\title{Deformation of finite dimensional $C^*$-quantum groupoids}
\author{Jean-Michel Vallin}
\address{Administrative address:  UMR CNRS 6628 Universit\'e D'Orl\'eans,
\vskip 0,3cm
  Postal address:  Institut de
Math\'ematiques de Chevaleret Projet alg\`ebres
\vskip 0,01cm
d'op\'erateurs et repr\'esentations, Plateau 7D, 175
rue du Chevaleret 75013 Paris
\vskip 0,3cm
  e-mail: jmva @ math.jussieu.fr}


\begin{abstract} In this work we prove, in a self contained way, that any finite dimensional $C^*$-quantum groupoid can be deformed in order that the square of the antipode is the  identity on the  base. We also prove that for any $C^*$-quantum groupoid with non abelian base, there is  uncountably many  $C^*$-quantum groupoids  with the same underlying algebra structure but which are not isomorphic to it. In fact, the $C^*$-quantum groupoids are closed in an analog of the procedure presented by  D.Nikshych ([N] 3.7) in a more general situation.
\end{abstract}
\maketitle
\newpage

\vskip 7cm
\def\al{\begin{align*} \end{align*}}

\newpage
\newenvironment{dm}{\hspace*{0,15in} {\it Proof:}}{\qed}
\section{ Flipper (separating) projections}
\label{intro}

\subsection{\bf{Notations}}
\label{notations}
In what follows,  $N$ is a finite dimensional von Neumann algebra, so 
$N$ is isomorphic to a sum of matrix algebras $\underset\gamma \oplus
M_{n_\gamma}$, we denote  the family of minimal central 
projections of $N$ by $\{ p^\gamma \}$,  we denote  a given family of matrix
units for $N$ by $\{e^\gamma_{i,j} / 1 \leq
i,j
\leq n_\gamma 
\}$. 

Let's denote by $tr_N$ the canonical (algebraic) trace on $N$ and $E_{Z(N)}$ the canonical conditional expectation $E: N \mapsto Z(N)$ (the center of $N$), so one has: $ tr_N \circ E_{Z(N)} = tr_N$, and $E_N$ is faithful. An easy calculation gives for any $n \in N: E_{Z(N)}(n) = \underset {\gamma,i,j }\sum
\frac{1}{n_\gamma}{e^\gamma_{i,j}}n e^\gamma_{j,i} = \underset {\gamma }\sum
\frac{1}{n_\gamma}tr_N(n p^\gamma) p^\gamma $

We shall denote  the  opposite  von Neumann algebra of
$N$ by $N^o$, so this is  $N$ with the opposite multiplication, hence a matrix
unit of
$N^o$ is given by the transposed  of $N's$ : $\{e^\gamma_{j,i} / 1 \leq i,j \leq
n_\gamma \}$. We shall not  in general distinguish $N$ and $N^o$, except when the  multiplication occurs.

\subsection{Definitions}
\label{cf }
{\it i) The application $m: N \otimes N \mapsto N$ is the multiplication map defined for any $x = \underset i \sum  x_i \otimes  y_i$ by $m(x) = \underset i \sum x_iy_i$.

ii) The elements $f$ in $N \otimes N$ such that  \underline { for the  multiplication of $N^o \otimes N$} , one has:  $f(a^o \otimes 1) = f(1 \otimes a)$ for any $a$ in $N$ and $m(f) =1$, are called the separating elements of $N$.}

\subsection{Example (cf [Val1],[Val2])}
\label{explosion}
The element $e = \underset \gamma \sum
\underset {i,j }\sum
\frac{1}{n_\gamma}{e^\gamma_{i,j}}^o\otimes e^\gamma_{j,i} $ is the only orthogonal 
projection  of $ N^o\otimes N$ such that, for any $n$ 
in $N$: $e(  n^o \otimes 1) = e(1 \otimes n)$ and   if
$e(1 \otimes n) = 0 $ then $n=0$. As an element of $N \otimes N$,  $e$ is also known as the symmetric separating element of $N$.

\subsection{Lemma}
\label{tortueux}
{\it For any  $y$ in $N \otimes N$,  \underline {and for the  multiplication of $N^o \otimes N$}, one has m(ey) = m(y).}
\newline
\begin{dm}
For any $x,z$ in $N$, one has:
\begin{align*}
m(e( x^o \otimes z))
&= m( \underset \gamma \sum
\underset {i,j }\sum
\frac{1}{n_\gamma}{(xe^\gamma_{i,j}})^o\otimes e^\gamma_{j,i}z) = \underset \gamma \sum
\underset {i,j }\sum
\frac{1}{n_\gamma}xe^\gamma_{i,j}e^\gamma_{j,i}z \\
&= x( \underset \gamma \sum
\underset {i,j }\sum
\frac{1}{n_\gamma}e^\gamma_{i,j}e^\gamma_{j,i})z 
= xz \\
&= m(x \otimes z)
\end{align*}
the lemma follows
\end{dm}

\subsection{Proposition}
\label{plaisant}
{\it Let $f$ be any element of $N \otimes N$, the following assertions are equivalent:

1) $f$ is a separating element of $N$,

2)$f$ is any {\underline{element}} of $N \otimes N$ such that in $N^o \otimes N$: $fe=f$ and $ef= e$ ,

3) $f$ is a (not necessarily orthogonal) projection of $N^o \otimes N$ with the same direction than $e$ 

4) their exists an element $g$ in $N$ such that $E_{Z(N)}(g)= 1$ and  $f = (1 \otimes g)e$ in $N^o \otimes N$ (or equivalently in $N \otimes N$).}
\newline
\begin{dm}
Let $f = \underset i \sum  f^o_i \otimes  f'_i$ and $h =  \underset i \sum  h^o_i \otimes  h'_i$ verifying 1) then in $N^o \otimes N$,  one has:
\begin{align*}
fh  
&= f(\underset i \sum  h^o_i \otimes  h'_i) =  \underset i \sum f(  h^o_i \otimes  h'_i) =  \underset i \sum f(  h^o_i \otimes 1) (1^o \otimes h'_i) \\
&=  \underset i \sum f( 1^o \otimes  h_i) (1^o \otimes h'_i) =    \underset i \sum f( 1^o \otimes  h_i) (1^o \otimes h'_i)  =  f( 1^o \otimes   \underset i \sum h_ih'_i) =  f 
\end{align*}
hence 1) implies 2).

Let $f = \underset i \sum  f^o_i \otimes  f'_i$  be any element verifying 2) then in $N^o \otimes N$:
\begin{align*}
f^2
&= fefe = f(ef)e = fe^2 = fe = f
\end{align*}
so $f$ is an algebraic projection, and as $fe= f$ (resp. $ef = e$), one has: $\ker e \subset \ker f$ (resp. $\ker f \subset \ker e$) and $f$ verifies 3).

Let  $f = \underset i \sum f^o_i \otimes  f'_i$  be any element  verifying 3), then one has:
\begin{align*}
f
&= fe =  \underset i \sum (  f^o_i \otimes  f'_i)e = \underset i \sum ( f^o_i \otimes  f'_i)e =  \underset i \sum (1^o \otimes f'_i)(  f^o_i \otimes 1)e\\
& =  \underset i \sum (1^o \otimes f'_i)( 1^o \otimes f_i)e =   (1^o \otimes \underset i \sum f'_if_i)e
\end{align*}

So there exists a $g$ in $N$ such that $f = (1^o \otimes g)e$ and one has:
\begin{align*}
e&= ef = e(1 \otimes g)e = (\underset \gamma \sum
\underset {i,j }\sum
\frac{1}{n_\gamma}{e^\gamma_{i,j}}^o\otimes e^\gamma_{j,i})(1 \otimes g)e \\
&= \underset \gamma \sum
\underset {i,j }\sum
(\frac{1}{n_\gamma}{e^\gamma_{i,j}}^o\otimes e^\gamma_{j,i}g)e \\
&=  \underset \gamma \sum
\underset {i,j }\sum
\frac{1}{n_\gamma}(1 \otimes e^\gamma_{j,i}g)({e^\gamma_{i,j}}^o\otimes 1)e = \underset \gamma \sum
\underset {i,j }\sum
\frac{1}{n_\gamma}(1 \otimes e^\gamma_{j,i}g)( 1 \otimes e^\gamma_{i,j})e \\
&= (1 \otimes \underset \gamma \sum
\underset {i,j }\sum
\frac{1}{n_\gamma}e^\gamma_{j,i}g e^\gamma_{i,j})e
\end{align*}

One deduces that $\underset \gamma \sum
\underset {i,j }\sum
\frac{1}{n_\gamma}e^\gamma_{j,i}g e^\gamma_{i,j} = 1$ and $f$ verifies 4).

To end let $f = \underset i \sum \overset o f_i \otimes  f'_i$  be any element  verifying 4) then for any $a$ in $N$, one has:
\begin{align*}
f(1 \otimes a) = (1 \otimes g)e(1 \otimes a) = (1 \otimes g)e(\overset o a \otimes 1) = f(\overset o a \otimes 1)
\end{align*}

Moreover: one easily checks that:  $m(f) = m((1 \otimes g)e) = \underset {i,j }\sum
\frac{1}{n_\gamma}e^\gamma_{j,i}g e^\gamma_{i,j} = 1$. So $f$ verifies 1).
\end{dm}

\subsection{\bf{Remarks}}
\label{courage}
What distinguishes $e$ among all separating elements of $N$ is the fact that it is an autoadjoint (orthogonal) projection for the natural $C^*$-structure of $N^o \otimes N$. When $N$ is not abelian there always exist  separating elements different from $e$: for example if $\gamma_0$ is such that $n_{\gamma_0} \not= 1$ then $f = ( 1 \otimes g)e$, where $g = 1 + e^{\gamma_0}_{1,1} - \frac{1}{n_{\gamma_0}-1}\underset {i \not= 1} \sum e^{\gamma_0}_{i,i}$, is such an element. The element $g$ in 4) needs not to be invertible, for example if $N = M_2(\mathbb C)$ (the $C^*$-algebra of complexe $2\times 2$ matrices) with it's usual matrix unit $(e_{i,j})$ then $g = 2 e_{1,1}$ works..

\subsection{Lemma and notations}
\label{approfondi}
{\it If $g$ is any strictly positive element of $N$, let's define a new involution denoted by $*_g$on the underlying algebra of $N$ by the formulae: $ x^{*_g} = gx^*g^{-1}$, then with this new $*$-involution one gives to this algebra a new structure of finite $C^*$-algebra we shall denote $N_g$.}
\newline
\begin{dm} Let''s prove that $N_g$ is actually a $C^*$-algebra. As any finite dimensional  $C^*$-algebra, N is isomorphic to a  finite sum of matrix algebras over $\mathbb C$, then one can suppose that : $N= M_n(\mathbb C)$ for some $n \in \mathbb N$ which  acts naturally on $H= \mathbb C^n$. Let's denote by $<,>$ the natural scalar product of $\mathbb C^n$, if one defines a new pairing  on $\mathbb C^n$, by the formula $<\xi,\eta>_g = <g^{-\frac{1}{2}}\xi,g^{-\frac{1}{2}}\eta>$, for any $\xi, \eta \in \mathbb C^n$, it's easy to see that it is a new scalar product on $\mathbb C^n$. Let $x$ be any element in $N$, then for any $\xi,\eta \in \mathbb C^n$, one has:
\begin{align*}
<x\xi,\eta>_g
&= <g^{-\frac{1}{2}}x\xi,g^{-\frac{1}{2}}\eta> = <x\xi,g^{-1}\eta> = < g^{-\frac{1}{2}}\xi, g^{-\frac{1}{2}} gx^*g^{-1}\eta> \\
&= < \xi, gx^*g^{-1}\eta>_g
\end{align*}
and so $x^{*_g} $ appears to be the adjoint of $x$ for this new scalar product and $N_g$ is a $C^*$-algebra.
\end{dm}

\subsection{Proposition}
\label{elabore}
{\it Let $g$ be any strictly positive element of $N$, such that $E_{Z(N)}(g) = 1$, and $f $any element of $N \otimes N$, then the following assertions are equivalent:

i) $f$ is the  separating element of $N$ associated to $g$ by proposition \ref{plaisant}

ii) $f$ is an orthogonal projection of $(N^o)_{(g^{\frac{1}{2}})^o} \otimes N_{g^{\frac{1}{2}}}$ with the same  direction  than $e$.}
\newline
\begin{dm} If i) is true by proposition \ref{plaisant}, $f$ is an idempotent of $N^o \otimes N$ with the same direction than $e$.  Let's prove that it is self adjoint for $(N^o)_{(g^{\frac{1}{2}})^o} \otimes N_{g^{\frac{1}{2}}}$, using the separability property for $e$, one has: 
\begin{align*}
f^{*_{g^{\frac{1}{2}}}}
&= ((1 \otimes g)e)^{*_{g^{\frac{1}{2}}}} = ((g^{\frac{1}{2}})^o \otimes g^{\frac{1}{2}})((1 \otimes g)e)^* ( (g^{-\frac{1}{2}})^o \otimes g^{-\frac{1}{2}}) \\
&= ((g^{\frac{1}{2}})^o \otimes g^{\frac{1}{2}})e(1 \otimes g)( (g^{-\frac{1}{2}})^o \otimes g^{-\frac{1}{2}}) = (1 \otimes g)e( (g^{-\frac{1}{2}})^o \otimes g^{\frac{1}{2}}) \\
&= (1 \otimes g)e = f,
\end{align*}
hence, ii) is true. Now let's suppose that $f$ verifies ii), then by proposition \ref{plaisant}, there exists a $k$ in $N$ such that $f= (1 \otimes k)e$, as $f$ is self adjoint, one has:
\begin{align*}
f
&= f^{*_{g^{\frac{1}{2}}}}
= ((1 \otimes k)e)^{*_{g^{\frac{1}{2}}}} 
= ((g^{\frac{1}{2}})^o \otimes g^{\frac{1}{2}})((1 \otimes k)e)^* ( (g^{-\frac{1}{2}})^o \otimes g^{-\frac{1}{2}}) \\
&= ((g^{\frac{1}{2}})^o \otimes g^{\frac{1}{2}})e(1 \otimes k^*)( (g^{-\frac{1}{2}})^o \otimes g^{-\frac{1}{2}}) = (1 \otimes g)e( (g^{-\frac{1}{2}})^o \otimes k^*g^{-\frac{1}{2}}) \\
&= (1 \otimes g)e( 1 \otimes g^{-\frac{1}{2}}k^*g^{-\frac{1}{2}}),
\end{align*}
so we can deduce that:
\begin{align*}
1
&= E_{Z(N)}(g) = m(f) = m((1 \otimes g)e( 1 \otimes g^{-\frac{1}{2}}k^*g^{-\frac{1}{2}})) \\
&= m((1 \otimes g)e)g^{-\frac{1}{2}}k^*g^{-\frac{1}{2}} = E_{Z(N)}(g)g^{-\frac{1}{2}}k^*g^{-\frac{1}{2}} \\
&= g^{-\frac{1}{2}}k^*g^{-\frac{1}{2}},
\end{align*}
hence $k= g^* = g$ and $f$ verifies i).
\end{dm}

\section{Deformation of $C^*$-quantum groupoids}
\label{decadix}

Let us recall the definition of a  $C^*$-quantum groupoid (or a weak Hopf 
$C^\star$ -algebra):

\subsection{\bf{Definition}}(G.B\"ohm, K.Szlach\'anyi, F.Nill) 
[BoSz], [BoSzNi]
\label{Bohm}

{\it A weak Hopf  $C^*$-algebra is a collection $(A, \Delta, \kappa,
\epsilon)$ where: $A$ is a finite-dimensional $C^*$-algebra, $\Delta: 
A \mapsto A \otimes A$ is a generalized coproduct, which means that: 
$(\Delta\otimes i)\Delta = (i \otimes \Delta)\Delta$, $\kappa$ is an
antipode on $A$, i.e., a linear application  from $A$ to $A$ such 
that $(\kappa \circ *)^2 = i$ (where $*$ is the involution on $A$), 
$\kappa(xy) = \kappa(y)\kappa(x)$ for every  $x,y$ in $A$ with
$(\kappa \otimes \kappa) \Delta = \varsigma\Delta \kappa$ (where 
$\varsigma$ is the usual flip on $A\otimes A$). 

We suppose also that 
$(m(\kappa \otimes i)\otimes i)(\Delta \otimes i)\Delta(x) = 
(1 \otimes x)\Delta(1)$ (where $m$ is the multiplication of
tensors, i.e., $m(a \otimes b) = ab$), and that $\epsilon$ is a   
counit, i.e., a positive linear form on $A$ such that
$(\epsilon \otimes i)\Delta = (i\otimes \epsilon)\Delta = i$, and 
for every  $x,y$ in $A$: $(\epsilon \otimes\epsilon) ((x \otimes 1) 
\Delta(1)(1 \otimes y)) = \epsilon(xy)$. }

\subsection{\bf{Results}} (L.Vainerman, D.Nikshych [NV1], [NV2],[BoSzNi])
\label{denuit}
{\it If $(A, \Delta, \kappa, \epsilon)$ is a weak Hopf  $C^*$-algebra,
then the following assertions are true:

0) The sets 
\[ A_t =\{ x\in A/\Delta(x) = \Delta(1)(x\otimes 1) = 
(x \otimes 1)\Delta(1) \}\] 
\[A_s= \{ x \in A / \Delta(x) = 
\Delta(1)(1\otimes x) = (1\otimes x)\Delta(1)\} \] 
are sub $C^*$-algebras of $A$; we call them  respectively target and 
source  Cartan subalgebra of $(A,\Delta)$.

1)  The application 
$\epsilon_t = m(i \otimes \kappa)\Delta$  takes values in $A_t$ 
and $\epsilon_s = 
m(\kappa\otimes i)\Delta$ takes values in $A_s$. We will call target
counit the application $\epsilon_t$ and we call source
counit $\epsilon_s$.

2) The $C^*$-algebra $A$ has a unique projection $p$, called the
Haar projection, characterized by the relations: $\kappa(p) = p$, 
$\epsilon_t(p) = 1$ and for every  $a$ in $A$, $ap =\epsilon_t(a)p$.

3) There exists a unique faithful positive linear form $\phi$, called 
the normalized Haar measure of $(A,\Delta,\kappa,\epsilon)$, satisfying
the following three properties:

$\phi\circ\kappa = \phi$, $(i \otimes \phi)(\Delta(1)) = 1$ and, for every 
$x,y$ in $A$:
\[(i\otimes\phi)((1\otimes y)\Delta(x)) = \kappa((i\otimes
\phi)(\Delta(y)(1\otimes x))).\]

4) $(\kappa \otimes i)\Delta(1)$ is a separating element for $N = A_t$,

5) One says that the collection $(A, \Delta, \kappa, \epsilon)$ is a weak 
Kac algebra if it is a $C^*$-quantum groupoid the antipode of  which
 is involutive, this is equivalent to the fact that $\phi$ is a
trace.}

\vskip 0.5cm

\subsection{\bf {Remark}}
{\it  The $C^*$-quantum groupoids we consider are not weak Kac algebras in
general, even there exist $C^*$-quantum groupoids for which $ 
\kappa^2 \not= Id$ on $A_t$ or $A_s$, but  one can deform the initial structure in order that $ 
\kappa^2(x) = x$, for any $x$ in $A_t$ or $A_s$. This is the object of what follows.}

\subsection{\bf{Proposition}}
\label{usure}
{\it If  $(A,\Delta,\kappa,\epsilon)$ is a $C^*$-quantum groupoid, there exists a unique {\bf invertible } element $q \in A_t$, such that   for any matrix unit   $( e^\gamma_{j,i})$ of $A_t$:
\vskip 0,1cm
1) $\Delta(1) = \underset \gamma \sum
\underset {i,j }\sum
\frac{1}{n_\gamma}{\kappa^{-1}(e^\gamma_{i,j}q)}\otimes e^\gamma_{j,i}$
\vskip 0,1cm
2) $q = \kappa^2(q)$ and $ E_{Z(A_t)}(q)= 1$
\vskip 0,1cm

3) For any $x \in A_tA_s$ (the $C^*$-algebra generated by $A_t$ and $A_s$), one has:  $\kappa^2(x) = q^{-1}\kappa(q)xq\kappa(q^{-1})$
\vskip 0,1cm
4) q is {\bf positive}.}
\vskip 0,2cm

\begin{dm} Let's use the notations of proposition \ref{plaisant} in the case of $N= A_t$: as $f=(\kappa \otimes i)\Delta(1)$ is a separatinq element for $N = A_t$, there exists a unique $q$ in $A_t$, such that $f = (1 \otimes q)e$ and for which 1) is true. Let's denote by $h$ any element in $A$ such that $hq = 0$, as  $\epsilon$ is a co-unity for $A$, one has:

\begin{align*}
h
&= h(\epsilon \otimes i)(\Delta(1)) =  (\epsilon \otimes i)((1 \otimes h)\Delta(1)) \\
&= (\epsilon \circ \kappa^{-1} \otimes i)(\kappa \otimes i)((1 \otimes h)\Delta(1)) = (\epsilon \circ \kappa^{-1} \otimes i)((1 \otimes h)(\kappa \otimes i)(\Delta(1))) \\
&= (\epsilon \circ \kappa^{-1} \otimes i)((1 \otimes hq)e)= 0
\end{align*}

Hence $q$ is invertible in $A$. From proposition \ref{plaisant}, 4) one has: $ \underset  \gamma \sum
\underset {i,j }\sum
\frac{1}{n_\gamma}{e^\gamma_{i,j}q}e^\gamma_{j,i} = 1$. 

As $\Delta(1)$ is selfadjoint one has:

\begin{align*}
\underset \gamma \sum
\underset {i,j }\sum
\frac{1}{n_\gamma}{\kappa^{-1}(e^\gamma_{i,j}q)}\otimes e^\gamma_{j,i}
&= \underset \gamma \sum
\underset {i,j }\sum
\frac{1}{n_\gamma}{\kappa^{-1}(e^\gamma_{i,j}q)^*}\otimes e^\gamma_{i,j}
\end{align*}

As the family $(e^\gamma_{i,j})$ is a base for $A_t$, one has for any $\gamma$ and $i,j$: 
$ \kappa^{-1}(e^\gamma_{i,j}q)
= \kappa^{-1}(e^\gamma_{j,i}q)^*
$
So for any $x$ in $A_t$: 
$$ \kappa^{-1}(xq)
= \kappa^{-1}(x^*q)^*,$$ applying this to  $x =1$ and $x=q$ one also has: $\kappa^{-1}(q)^* = \kappa^{-1}(q)$, and $\kappa^{-1}(q^2) = \kappa^{-1}(q^*q)^*$. so  using the fact that $( \kappa \circ \star)^2 =1$, one obtains:
\begin{align*}
\kappa^{-1}(q)
&= \kappa^{-1}(q^2)\kappa^{-1}(q^{-1}) = \kappa^{-1}(q^*q)^*\kappa^{-1}(q^{-1}) \\
&=  (\kappa^{-1}(q)\kappa^{-1}(q^*))^*\kappa^{-1}(q^{-1}) \\
&=  \kappa^{-1}(q^*)^*\kappa^{-1}(q)^*\kappa^{-1}(q^{-1}) = \kappa^{-1}(q^*)^*\kappa^{-1}(q)\kappa^{-1}(q^{-1}) = \kappa^{-1}(q^*)^* \\
&= \kappa(q),
\end{align*}
hence      
$\kappa^2(q)= q$ and $q$ verifies 2). 

The fact that for any $x$ in $A_t$ one has $\kappa^{-1}(xq)
= \kappa^{-1}(x^*q)^*$  implies that:
\begin{align*}
\kappa^2(x) 
&= \kappa(\kappa(x)) = \kappa({\kappa}^{-1}(x^*)^*) = \kappa({\kappa}^{-1}(x^*)^*\kappa^{-1}(q)  \kappa^{-1}(q^{-1})) \\
&= q^{-1}  \kappa({\kappa}^{-1}(x^*)^*\kappa^{-1}(q) )= q^{-1} \kappa({\kappa}^{-1}(x^*)^* \kappa^{-1}(q)^*)\\
&= q^{-1}\kappa((\kappa^{-1}(q){\kappa}^{-1}(x^*))^* ) = q^{-1}\kappa({\kappa}^{-1}(x^*q)^* ) \\
&= q^{-1} \kappa({\kappa}^{-1}(xq) ) \\
&= q^{-1}xq
\end{align*}

Now let $y$ be any element of $A_s$, and let $z$ be the element in $A_t$ such that $y = \kappa(z)$, then, as $A_s$, and  $A_t$ commute, for any $x$ in $A_t$, one has:
\begin{align*}
\kappa^2(xy)
&= \kappa^2(y)\kappa^2(x) = \kappa(\kappa^2(z))\kappa^2(x) =
 \kappa(q^{-1}\kappa(z)q)q^{-1}xq  \\
&= \kappa(q)\kappa(z)\kappa(q^{-1})q^{-1}xq = \kappa(q)q^{-1}yxq\kappa(q^{-1}) \\
&= \kappa(q)q^{-1}xyq\kappa(q^{-1}),
\end{align*}
so $q$ verifies 3)

By an other calculus, one has:
\begin{align*}
\kappa^{-1}(q)
&=\kappa^{-1}(q)^* = (\kappa^{-1}(q)^*)^{-1}\kappa^{-1}(q)^*\kappa^{-1}(q)^* = (\kappa^{-1}(q))^{-1}\kappa^{-1}(q^{2})^* \\
&= \kappa^{-1}(q)^{-1}\kappa^{-1}(q)\kappa^{-1}(q^*) \\
&= \kappa^{-1}(q^*),
\end{align*}

So $q = q^*$, if $|q |$ denotes the module of $q$ in his polar decomposition, as a consequence $q$ and $|q |^{\frac{1}{2}}$ commute, hence $\kappa^2(|q|^{\frac{1}{2}}) = |q|^{\frac{1}{2}}$. But one has: $(\kappa \circ *)^2 = 1$, so $\kappa((\kappa(|q|^{\frac{1}{2}})^*) = |q|^{\frac{1}{2}} = \kappa((\kappa(|q|^{\frac{1}{2}}))$, this implies that $\kappa(|q|^{\frac{1}{2}})$ is hermitian. Now one has:  $\kappa(|q|) = \kappa(|q|^{\frac{1}{2}})\kappa(|q|^{\frac{1}{2}} )= \kappa(|q|^{\frac{1}{2}})^*\kappa(|q|^{\frac{1}{2}})  \geq 0$.

One can view $A_s$ as a sub-$C^*$-algebra of the linear operators of  a finite Hilbert space $H$, and $A_t$ as acting on $\underset \gamma \oplus \mathbb C^{n_\gamma}$ in such a way that the $e_{i,j}^\gamma$'s form the canonical matrix unit of $\mathcal L(\mathbb C^{n_\gamma})$. Let's denote by $\epsilon_{i}^\gamma$ the canonical base of $\mathbb C^{n_\gamma}$. As  $\Delta(1)$ is positive, then for any $\eta$ in $H$ and any $\epsilon_{k}^{\gamma^0}$, one has:
\begin{align*}
0
& \leq (\Delta(1)(\eta \otimes \epsilon_{k}^{\gamma^0}), \eta \otimes \epsilon_{k}^{\gamma^0})\\
& \leq \underset \gamma \sum
\underset {i,j }\sum
\frac{1}{n_\gamma} ({\kappa^{-1}(e^\gamma_{i,j}q)\eta}\otimes e^\gamma_{j,i}\epsilon_{k}^{\gamma^0}, \eta \otimes \epsilon_{k}^{\gamma^0}) \\
& \leq \frac{1}{n_{\gamma^0}}
\underset {j }\sum
 ({\kappa^{-1}(e^{\gamma^0}_{k,j}q)\eta}\otimes \epsilon_{j}^{\gamma^0}, \eta \otimes \epsilon_{k}^{\gamma^0}) \\
& \leq \frac{1}{n_{\gamma^0}}
\underset {j }\sum
 (\kappa^{-1}(e^{\gamma^0}_{k,j}q)\eta, \eta)(  \epsilon_{j}^{\gamma^0},\epsilon_{k}^{\gamma^0})\\
& \leq \frac{1}{n_{\gamma^0}}
 (\kappa^{-1}(e^{\gamma^0}_{k,k}q)\eta, \eta),
\end{align*} 

One deduces that   $\kappa^{-1}(e^{\gamma^0}_{k,k}q) \geq 0$, and $ \kappa(q) = \kappa^{-1}(q) = \underset {\gamma^0} \sum \underset {k} \sum \kappa^{-1}(e^{\gamma^0}_{k,k}q) \geq 0$.  As $q$ and $|q|$ commute, $\kappa(q)$ and $\kappa(|q|)$ are commuting each other, positive and also:
$$\kappa(q)^2 = \kappa(q^2) = \kappa(q^*q) = \kappa(|q|^2) = \kappa(|q|)^2,$$
hence one has  $\kappa(q) = \kappa(|q|)$, so $q = |q| \geq 0$ and $q$ verifies 4). \end{dm}

\subsection{Lemma}
\label{labeur}
{\it   Let $q$ be the canonical element of $A_t$ defined in proposition \ref{usure},. If the base $A_t$ is not commutative then the  set of  strictly positive elements $k$ in $A_t$, such that $\kappa^2(k) = k$ and $E_{Z(A_t)}(k^{-1}q) = 1$ is uncountable, more precisely the set of the $k^{-1}q$'s spectra is uncountable}
\newline
\begin{dm}
Due to proposition \ref{usure} 3), a strictly positive  element $k$ of $A_t$ such that $\kappa^2(k) = k$ is just a strictly positive element of  commuting with $q$, let's use the notations of \ref{notations} for $N= A_t$, for any $\gamma$, $q^\gamma $ and $k^\gamma $ appear to be two strictly positive matrices in $M_{n_\gamma}(\mathbb C)$ which commute, up to conjugacy, these are two diagonal matrices with all diagonal  elements strictly positive, let's say $(k^\gamma_i)_{i= 1..n}$ and $(q^\gamma_i)_{i= 1..n}$. The equality $E_{Z(N)}(k^{-1}q) = 1$ just means that $\underset i \sum ( k_i^\gamma)^{-1}q_i^{\gamma} = n_\gamma$. The lemma follows. \end{dm}

\subsection{\bf{Theorem}}
\label{deformation}
{\it  Let $(A,\Delta, \kappa, \epsilon)$ be any $C^*$-quantum groupoid, and let $k$ be any positive invertible element of $A_t$ given by lemma \ref{labeur}. Let's give the algebra $A$ a new involution $*{_k}$, a new coproduct $\Delta_k$, a new antipode $\kappa_k$ and a new counit $\epsilon_k$ by the following formulas: 

$\forall a \in A$,  $a^{*_k} =  ka^*k^{-1}$,  $\Delta_k(a) = \Delta(a)(1\otimes k^{-1})$, $\kappa_k(a) = k\kappa(x)k^{-1}$, and $\epsilon_k(a) = \epsilon(ak).$

If one denotes by $A_k$, the algebra $A$ together with the new involution $^{*_k}$, then the $4$-tuple $(A_k,\Delta_k, \kappa_k, \epsilon_k)$ is a new $C^*$-quantum groupoid, with the same underlying Cartan subalgebras than $(A,\Delta, \kappa, \epsilon)$ and  $k^{-1}q$ is the canonical element associated by proposition \ref{usure} to $(A_k,\Delta_k, \kappa_k, \epsilon_k)$  .}
\newline
\begin{dm}  By lemma \ref{approfondi}, $A_k$ is a $C^*$-algebra. For any $a,b \in A$,  using the notations of proposition  \ref{usure} and the fact that $\kappa^{-1} \otimes i$ is multiplicative viewed as an application from $N^o \otimes N$ to $N \otimes N$, one has: 
\begin{align*}
\Delta_k(a)\Delta_k(b) 
&= \Delta(a)(\Delta(1)(1 \otimes k^{-1})\Delta(1))\Delta(b)(1 \otimes k^{-1}) 
\\
&= \Delta(a)(\kappa^{-1} \otimes i)(f(1 \otimes k^{-1})f)\Delta(b)(1 \otimes k^{-1})  \\
&= \Delta(a)(\kappa^{-1} \otimes i)(f(1 \otimes k^{-1}q)e)\Delta(b)(1 \otimes k^{-1}) \\
&= \Delta(a)(\kappa^{-1} \otimes i)(f(1 \otimes E_{Z(N)}(k^{-1}q))\Delta(b)(1 \otimes k^{-1}) \\
&= \Delta(a)(\kappa^{-1} \otimes i)(f)\Delta(b)(1 \otimes k^{-1}) =  \Delta(a)\Delta(1)\Delta(b)(1 \otimes k^{-1}) \\
&= \Delta_k(ab),
\end{align*}
so $\Delta_k$ is multiplicative and obviously linear. As $A_s$ and $A_t$ are commuting, then one also has:
\begin{align*}
\Delta_k(a^{*_k})
&= \Delta_k(ka^*k^{-1})=  \Delta(ka^*k^{-1})(1 \otimes k^{-1}) =  \Delta(k)\Delta(a^*)\Delta(k^{-1})(1 \otimes k^{-1})\\
&= (k \otimes 1)\Delta(a^*)(k^{-1}   \otimes k^{-1}) =  (k \otimes k)(1\otimes k^{-1})\Delta(a)^*(k^{-1}   \otimes k^{-1}) \\
&= (k \otimes k)(\Delta(a)(1\otimes k^{-1}))^*(k^{-1}   \otimes k^{-1}) = (k \otimes k)\Delta_k(a)^*(k^{-1}   \otimes k^{-1}) \\
&= \Delta_k(a)^{*_k},
\end{align*}
hence  $\Delta_k$ is a $*$-morphism. In an other hand, we have:

\begin{align*}
(\Delta_k \otimes i)\Delta_k(a) 
&= (\Delta_k \otimes i)(\Delta(a)(1 \otimes k^{-1})) = (\Delta \otimes i)(\Delta(a))(1 \otimes  k^{-1} \otimes k^{-1}) \\
&= (i \otimes \Delta)(\Delta_k(a)(1 \otimes  k))(1 \otimes  k^{-1} \otimes k^{-1}) \\
&= (i \otimes \Delta)(\Delta_k(a))(1 \otimes \Delta(1)) (1 \otimes k \otimes  1)(1 \otimes  k^{-1} \otimes k^{-1})  \\
&= (i \otimes \Delta)(\Delta_k(a))(i \otimes \Delta)(1\otimes 1) (1 \otimes  1 \otimes k^{-1})  \\
&= (i \otimes \Delta)(\Delta_k(a)) (1 \otimes  1\otimes k^{-1})  \\
&= (i \otimes \Delta_k)(\Delta_k(a)) 
\end{align*}

So $\Delta_k$ is a coproduct. As $A_t$ and $A_s$ commute, one has:
\begin{align*}
\varsigma(\kappa_k  \otimes \kappa_k )\Delta_k(a)
&= (k \otimes k)\varsigma(\kappa  \otimes \kappa )\Delta_k(a)(k^{-1} \otimes k^{-1})  \\
&= (k \otimes k)\varsigma((\kappa  \otimes \kappa )(\Delta(a) (1\otimes k^{-1}))(k^{-1} \otimes k^{-1}) \\
&=  (k\otimes  k)(\kappa(k^{-1})\otimes 1)\varsigma((\kappa  \otimes \kappa )\Delta(a))(k^{-1} \otimes k^{-1}) \\
&= (k\otimes  k)(\kappa(k^{-1})\otimes 1)\Delta(\kappa(a))(k^{-1} \otimes k^{-1}) \\
&= (\kappa(k^{-1})k\otimes  k))\Delta(\kappa(a))(k^{-1}\otimes k^{-1}) \\
&= (\kappa(k^{-1}) \otimes k)\Delta(k\kappa(a)k^{-1})(1\otimes  k^{-1}) \\
&=  (\kappa(k^{-1}) \otimes k)\Delta(\kappa_k(a))(1\otimes  k^{-1})\\
&=  (\kappa(k^{-1}) \otimes k)\Delta(1)\Delta(\kappa_k(a))(1\otimes  k^{-1})
\end{align*}

But, as  $\kappa^{-1} \otimes i$ is multiplicative from  $A^o \otimes A$ to $A \otimes A$,  structure, and using the fact that $\kappa^2(k) = k$, one has:
\begin{align*}
(\kappa(k^{-1}) \otimes k)\Delta(1)
&= (\kappa^{-1}(k^{-1}) \otimes k)\Delta(1) =  (\kappa^{-1} \otimes i)(( \overset o  {k}^{-1} \otimes k)(\kappa \otimes i)(\Delta(1)) \\
&= (\kappa^{-1} \otimes i)(\kappa \otimes i)(\Delta(1)) = \Delta(1)
\end{align*}

Replacing this in the former list of equalities, one obtains:
\begin{align*}
\varsigma(\kappa_k  \otimes \kappa_k )\Delta_k(a)
&= \Delta(1)\Delta(\kappa_k(a))(1\otimes  k^{-1}) = \Delta_k(\kappa_k(a))
\end{align*}

Hence $\kappa$ is an antipode. For every $a$ in $A$, one has:
\begin{align*}
\kappa_k((\kappa_k(a^{*_k})^{*_k})
&= k\kappa (\kappa_k(a^{*_k})^{*_k})k^{-1} = k\kappa (k(\kappa_k(a^{*_k})^*k^{-1})k^{-1} \\
&= k\kappa(k^{-1})\kappa (\kappa_k(a^{*_k})^*)\kappa(k)k^{-1} \\
&= k \kappa(k^{-1})\kappa ((k\kappa (a^{*_k})k^{-1})^*)\kappa(k)k^{-1} \\
&= k \kappa(k^{-1})\kappa (k^{-1}\kappa (a^{*_k})^*k)\kappa(k)k^{-1} =  k \kappa (\kappa (a^{*_k})^*)k^{-1} \\
&= k \kappa (\kappa (ka^{*}k^{-1})^*)k^{-1} =  k \kappa (\kappa(k)^*\kappa (a^{*})^*\kappa(k^{-1})^*)k^{-1} \\
&=   k \kappa(\kappa(k^{-1})^*) \kappa(\kappa (a^{*})^*)\kappa(\kappa(k)^*)k^{-1} \\
&=   k \kappa(\kappa((k^{-1})^*)^*) a\kappa(\kappa(k^*)^*)k^{-1} =  k k^{-1}akk^{-1} \\
&= a.
\end{align*}
Let's now use the co-unity  $\epsilon$  and  the fact that $k$ commutes with $A_s$:
\begin{align*}
(\epsilon_k \otimes i)\Delta_k(a) 
&= (\epsilon \otimes i)(  \Delta(a)(k \otimes k^{-1})) = (\epsilon \otimes i)  \Delta(ak)k^{-1} = akk^{-1} = a
\end{align*}
and more obviously:
\begin{align*}
(i \otimes \epsilon_k)\Delta_k(a) 
&= (i \otimes \epsilon)(  \Delta(a)(1 \otimes  k^{-1}k)) = (i \otimes \epsilon ) \Delta(a)= a,
\end{align*}
hence $\epsilon_k$ is also a co-unity, and for every $a,b$ in $A$, one has:
\begin{align*}
\epsilon_k(xy) 
&= \epsilon(xyk) = \epsilon(xkk^{-1}yk) =(\epsilon \otimes \epsilon)((x \otimes 1)\Delta(1)( k \otimes k^{-1}yk)) \\
&= (\epsilon_k \otimes \epsilon_k)((x \otimes 1)\Delta(1)( 1 \otimes k^{-1}y)) \\
&= (\epsilon_k \otimes \epsilon_k)((x \otimes 1)\Delta_k(1)( 1 \otimes y))
\end{align*} 

Then       for every $a$ in $A$, one has:
\begin{align*}
\epsilon_k(a^{*_k}a)
&= \epsilon(a^{*_k}ak)=  \epsilon(ka^{*}k^{-1}ak) = \epsilon(ka^{*}k^{-\frac{1}{2}}k^{-\frac{1}{2}}ak) \\
&= \epsilon((k^{-\frac{1}{2}}ak)^*(k^{-\frac{1}{2}}ak))\\
&\geq 0.
\end{align*}

Finally, using lemma \ref{tortueux}, for any $x$ in $A$, one has:
\begin{align*}
(m(\kappa_k &\otimes i) \otimes i)(\Delta_k \otimes i)\Delta_k(x)
\\
&= (m(\kappa_k \otimes i) \otimes i)((\Delta \otimes i)\Delta(x)( 1 \otimes k^{-1} \otimes k^{-1})) \\
&= (m \otimes i)((k \otimes 1 \otimes 1)(\kappa \otimes i \otimes i)((\Delta \otimes i)\Delta(x))( k^{-1}  \otimes k^{-1} \otimes k^{-1})) \\
&= (k \otimes 1)(m \otimes i)((\kappa  \otimes i \otimes i)((\Delta \otimes i)\Delta(x))( k^{-1}  \otimes 1 \otimes 1))(k^{-1} \otimes k^{-1})
\end{align*}

but $\kappa \otimes i \otimes i$ is multiplicative from $A\otimes A \otimes A$ to $A^o \otimes A \otimes A $, so we can write, using lemma \ref{tortueux}, that:
\begin{align*}
(m &\otimes i)((\kappa  \otimes i \otimes i)
((\Delta \otimes i)\Delta(x))( k^{-1}  \otimes 1\otimes 1)) \\
&= (m \otimes i)((\kappa  \otimes i \otimes i)
((\Delta \otimes i)\Delta(x)(\Delta(1) \otimes 1))( k^{-1}  \otimes 1\otimes 1))\\
&= (m \otimes i)(( \overset o k^{-1} \otimes 1 \otimes 1)(\kappa  \otimes i \otimes i)(\Delta(1) \otimes 1)(\kappa  \otimes i \otimes i)
((\Delta \otimes i)\Delta(x)))\\
&= (m \otimes i)(( \overset o k^{-1} \otimes 1 \otimes 1)( f \otimes 1)(\kappa  \otimes i \otimes i)
((\Delta \otimes i)\Delta(x))) \\
&= (m \otimes i)(( 1 \otimes  k^{-1} \otimes 1)( f \otimes 1)(\kappa  \otimes i \otimes i)
((\Delta \otimes i)\Delta(x)))\\
&= (m \otimes i)(( e \otimes 1)(\kappa  \otimes i \otimes i)
((\Delta \otimes i)\Delta(x))) \\
&= (m \otimes i)((\kappa  \otimes i \otimes i)
((\Delta \otimes i)\Delta(x))) \\
&= (1 \otimes x)\Delta(1),
\end{align*}

Hence replacing this in the former equality and using the fact that $k$ commutes with $A_s$,
one has:

\begin{align*}
(m(\kappa_k &\otimes i) \otimes i)(\Delta_k \otimes i)\Delta_k(x)\\
&= (k \otimes 1)(1 \otimes x)\Delta(1)(k^{-1} \otimes k^{-1})\\
&= (1 \otimes x)\Delta(1)(1 \otimes k^{-1}) \\
&= (1 \otimes x)\Delta_k(1)
\end{align*}

So $(A_k,\Delta_k, \kappa_k, \epsilon_k)$ is a $C^*$-quantum groupoid. As for any $x$ in $A$, one has $\Delta_k(x) = \Delta(x)(1 \otimes k^{-1})$, then it's obvious that the Cartan subalgebras are the same with this new structure. Finally there is a single element of $A_t$ verifying proposition \ref{usure} 1), but, as $\Delta_k(1) = \Delta(1)(1 \otimes k^{-1})$, an  easy computation gives that $k^{-1}q$ verifies proposition \ref{usure} for $\Delta_k$, the theorem follows.
\end{dm} 

\subsection{Corollary}
\label{butessentiel}
{\it  Every finite dimensional $C^*$-quantum groupoid  can be deformed in such a way that the antipode becomes involutive on the Cartan subalgebras. }
\newline
\begin{dm}
It suffices to apply the last theorem for $k= q$.
\end{dm}

\subsection{Corollary}
\label{drole}
{\it  Any  $C^*$-quantum groupoid,  the Cartan subalgebras of which are not abelian, can be deformed in such a way that the antipode becomes {\bf non}-involutive on the Cartan subalgebras;  there  is  uncountably many  non isomorphic $C^*$-quantum groupoids with its underlying algebra structure. }
\newline
\begin{dm}
Let $(A,\Delta, \kappa, \epsilon)$ be any  $C^*$-quantum groupoid,  for which $A_t$ is  not abelian, then obviously there exits a $k'$ in $A_t$ given by lemma \ref{labeur} such that $k'^{-1}q$ is not in the center of $A_t$, then by theorem \ref{deformation} the antipode $\kappa_{k'}$ is not involutive on $A_t$.

As the canonical element associated with $(A_k,\Delta_k, \kappa_k, \epsilon_k)$ by proposition \ref{usure} is $k^{-1}q$, the spectrum of  $k^{-1}q$ is an invariant of the isomorphism class of $(A_k,\Delta_k, \kappa_k, \epsilon_k)$, one can conclude using lemma \ref{labeur}.
\end{dm}

\subsection{Remark}
{\it  As the deformation moves the natural involution of $A$, this leads to a natural question: for any  $C^*$-quantum groupoid,  the Cartan subalgebras of which are not abelian,  is  there  uncountably many  non isomorphic $C^*$-quantum groupoids with its underlying $C^*$-algebra structure?}

\vskip 2cm

REFERENCES

[BoSz] G. B\"OHM \& K.SZLACH\'ANYI, Weak C*-Hopf algebras: the  
coassociative symmetry of non integral dimensions, {\it  Quantum groups 
and quantum spaces. Banach Center Publications} {\bf 40} (1997), 9-19.

[BoSzNi] G. B\"OHM, K.SZLACH\'ANYI \& F.NILL, Weak Hopf Algebras I. 
Integral Theory and $C^*$-structure. {\it Journal of Algebra} {\bf 221} 
(1999),  385-438.

[N] NIKSHYCH D On the structure of weak Hopf algebras, {\it
Advances in Mathematics} {\bf 170} (2002), 257-286.

[NV1] D. NIKSHYCH \& L. VAINERMAN, Algebraic versions of a 
finite-dimensional quantum groupoid. {\it Lecture Notes in Pure and 
Appl. Math.} {\bf 209} (2000), 189-221.

[NV2] D. NIKSHYCH \& L. VAINERMAN, Finite Quantum Groupoids and Their
Applications, in New Directions in Hopf Algebras {\it MSRI Publications} {\bf 43} (2002) Cambridge University Press , 211-262.

[Val1] J.M. VALLIN, Groupo\" \i des quantiques finis. {\it Journal 
of Algebra} {\bf 26} (2001), 425-488.

[Val2]  J.M.  VALLIN  Multiplicative partial
isometries and finite quantum groupoids : Proceedings of the
Meeting of Theoretical Physicists and Mathematicians,
Strasbourg, 2002. IRMA Lectures in Mathematics and
Theoritical Physics {\bf 2} 189-227. 

\end{document}